\documentclass[11pt,reqno]{amsart}
\usepackage{amsmath,amssymb,amsthm,mathrsfs,dsfont}
\usepackage{CJK}
\usepackage{amsmath}
\usepackage{dsfont}
\usepackage{mathrsfs}
\usepackage{amsmath,amssymb}
\usepackage{amsfonts}
\usepackage{hyperref}
\usepackage{amsthm}
\usepackage{graphicx}
\usepackage{subfigure}
\usepackage{xcolor}
\usepackage{overpic}

\newtheorem{theorem}{Theorem}[section]
\newtheorem{lemma}{Lemma}[section]
\newtheorem{definition}{Definition}[section]
\newtheorem{proposition}{Proposition}[section]
\numberwithin{equation}{section}
\newtheorem{remark}{Remark}[section]

\def\pf{{\textit {Proof:} }}

\newfont{\bb}{msbm10 at 12pt}

\newcommand{\bal}{\begin{aligned}}      \newcommand{\eal}{\end{aligned}}
\newcommand{\ba}{\begin{array}}      \newcommand{\ea}{\end{array}}
\newcommand{\bc}{\begin{center}}     \newcommand{\ec}{\end{center}}
\newcommand{\be}{\begin{enumerate}}  \newcommand{\ee}{\end{enumerate}}
\newcommand{\beq}{\begin{eqnarray}}  \newcommand{\eeq}{\end{eqnarray}}
\newcommand{\beQ}{\begin{eqnarray*}} \newcommand{\eeQ}{\end{eqnarray*}}
\newcommand{\bi}{\begin{itemize}}    \newcommand{\ei}{\end{itemize}}
\newcommand{\bt}{\begin{tabular}}    \newcommand{\et}{\end{tabular}}
\newcommand{\bdm}{\begin{displaymath}} \newcommand{\edm}{\end{displaymath}}

\def\qed{\hfill{Q.E.D.}\smallskip}

\title{Combinatorial Calabi flow with surgery on surfaces}
\author{Xiang Zhu, Xu Xu}
\address[Xiang Zhu]{Yau Mathematical Sciences Center, Tsinghua University, Beijing 100084, P.R. China}
\email{zhu-x13@mails.tsinghua.edu.cn}

\address[Xu Xu]{School of Mathematics and Statistics, Wuhan University, Wuhan 430072, P.R. China}
\email{xuxu2@whu.edu.cn}
\date{}

\begin{document}
\maketitle

\begin{abstract}
Computing uniformization maps for surfaces has been a challenging problem and has many practical applications.
In this paper, we provide a theoretically rigorous algorithm to compute such maps via
combinatorial Calabi flow for vertex scaling of polyhedral metrics on surfaces,
which is an analogue of the combinatorial Yamabe flow introduced by Luo \cite{L1}.
To handle the singularies along the combinatorial Calabi flow, we do surgery on the flow by flipping.
Using the discrete conformal theory established in \cite{GGLSW,GLSW},
we prove that for any initial Euclidean or hyperbolic polyhedral metric on a closed surface,
the combinatorial Calabi flow with surgery exists for all time and
converges exponentially fast after finite number of surgeries. The convergence is
independent of the combinatorial structure of the initial
triangulation on the surface.
\end{abstract}



\section{Introduction}\label{section 1}
\subsection{Backgrounds and main results}
One of the central topics in modern geometry concerns with
the canonical metrics on a given manifold, which is related to
special geometric structures on manifolds.
The flow method is an important approach for such problems.
To achieve this goal, Hamilton \cite{H1} introduced the Ricci flow and
Calabi \cite{Ca, Ca2} introduced the Calabi flow.
The Ricci flow has been
used to solve many longstanding problems in differential geometry, including the
Poincar\'{e} conjecture \cite{P1, P2, P3} and the sphere theorem \cite{BrSc}.
The Calabi flow is
\begin{equation}\label{Calabi flow}
\frac{\partial g}{\partial t}=\Delta_g K\cdot g
\end{equation}
on a Riemannian surface, where $g$ is the Riemannian metric, $K$ is the Gauss curvature and
$\Delta_g$ is the Laplace-Beltrami operator of $g$.
The long time existence and convergence of the Calabi flow (\ref{Calabi flow}) on closed Riemannian surfaces
were established in \cite{CXX, CP}.
The readers can also refer to \cite{CSC,CSC2,CLT} and the references therein for more information on Calabi flow on surfaces.

%Chru\'{s}ciel \cite{CP} proved that the Calabi flow (\ref{Calabi flow}) exists
%for all time and
%converges to a constant Gauss curvature metric on closed Riemannian surfaces using the Bondi mass estimation,
%assuming the existence of the constant
%Gauss curvature which is ensured by the uniformization theorem.
%Chen \cite{CXX} gave a geometrical proof of the long time existence and the convergence of the Calabi flow on closed surfaces.

However, computing special geometric structures on manifolds has been a challenging problem, even for the uniformization maps on surfaces.
Different from the smooth approach, one works on polyhedral manifolds to compute geometric structures on manifolds.
This idea dates at least back to Thurston \cite{T1},
who used the circle packing metric (a type of polyhedral metric) to study hyperbolic structures on three dimensional manifolds.
Thurston (\cite{T1}, Section 13.7) obtained the combinatorial obstruction for
the existence of constant combinatorial curvature circle packing metrics.
Motivated by Ricci flow on surfaces \cite{CH1, H2}, Chow-Luo \cite{CL} introduced the combinatorial surface
Ricci flow for Thurston's circle packing metrics on triangulated surfaces, which provides an effective way to
compute metrics with prescribed combinatorial curvatures (specially the uniformization maps) and has
many applications (see for example \cite{ZG} and the references therein).
Inspire by \cite{CL},
Ge \cite{Ge-thesis, Ge} introduced a combinatorial Calabi flow for Thurston's Euclidean circle packing metrics on triangulated surfaces
and proved the longtime existence and convergence of the flow.
Then Ge and the second author \cite{GX2} and Ge-Hua \cite{GH} studied the convergence of combinatorial Calabi flow for Thurston's hyperbolic circle packing metrics
on triangulated surfaces.
Similar to combinatorial Ricci flow,
the combinatorial Calabi flow can also be used to compute circle packing metrics with prescribed combinatorial curvatures \cite{Zhao}.
Combinatorial Ricci flow and Calabi flow are further used to
study $\alpha$-curvatures of circle packings on triangulated surfaces and sphere packings on 3-dimensional manifolds \cite{GX1, GXcv, GX3, GX4}.

To study the discrete conformal geometry of polyhedral metrics on maniflods,
R\v{o}cek-Williams \cite{RW} and Luo \cite{L1} introduced a notion
of discrete conformality for Euclidean polyhedral metrics (piecewise linear metrics or PL metrics for short) on triangulated surfaces independently,
which is now called vertex scaling.
Luo \cite{L1} further defined the combinatorial Yamabe flow for PL metrics.
%Bobenko-Pinkall-Springborn \cite{BPS} studied the vertex scaling introduced by Luo and R\v{o}cek-Williams
%and proved the global rigidity of combinatorial curvature with respect to the conformal factor, which was conjectured by Luo \cite{L1}.
%They further introduced vertex scaling for piecewise hyperbolic metrics and piecewise spherical metrics on triangulated surfaces.
Based on Bobenko-Pinkall-Springborn's work \cite{BPS} on vertex scaling and Penner's work \cite{Penner} on decorated Teichim\"{u}ller spaces,
Gu-Luo-Sun-Wu \cite{GLSW} and Gu-Guo-Luo-Sun-Wu \cite{GGLSW} recently proved a discrete uniformization theorem for Euclidean and hyperbolic polyhedral metrics on surfaces respectively,
which provides a constructive proof of the classical uniformization theorem on closed surfaces.
Combinatorial Yamabe flows with surgery were introduced in \cite{GGLSW, GLSW},
where the long-time existence and convergence were proved.
The finiteness of surgeries along the combinatorial Yamabe flow was proved by Wu \cite{Wu}.
The combinatorial Yamabe flow with surgery provides an effective algorithm to compute the uniformization maps on surfaces \cite{SWGL}.
Unlike the circle packing case \cite{CL,Ge-thesis,Ge,GH,GX2}, the convergence of the combinatorial Yamabe flow with surgery in \cite{GGLSW, GLSW}
is independent of the combinatorial structure of the initial triangulation on the surface.
Following \cite{L1}, Ge \cite{Ge-thesis} introduced the combinatorial Calabi flow for vertex scaling of PL metrics on triangulated surfaces
and proved the short time existence.
In this paper, by doing surgery along the combinatorial Calabi flow by flipping,
we prove the combinatorial Calabi flow with surgery for vertex scaling exists for all time and converges as time tends to infinity.
One of the main results of the paper is the following theorem on Euclidean combinatorial Calabi flow with surgery.
\begin{theorem}\label{main theorem for Calabi flow with surgery intro}
Suppose $(S, V)$ is a closed connected marked surface and $d_0$ is any initial piecewise linear metric on $(S, V)$.
Then the combinatorial Calabi flow with surgery
exists for all time and
converges exponentially fast to a piecewise linear metric $d^*$ with
constant combinatorial curvature after finite number of surgeries.
\end{theorem}

Note that the convergence is independent of the combinatorial structure of the triangulation on the surface,
which improves the convergence of combinatorial curvature flows for circle packing metrics on surfaces \cite{CL,Ge-thesis,Ge,GH,GX2}.
The combinatorial Calabi flow with surgery gives an effective way to compute polyhedral metrics on a surface with given combinatorial curvatures \cite{SLZXG}.
Especially, the combinatorial Calabi flow provides an algorithm to compute the uniformization maps on surfaces.

We further define the hyperbolic combinatorial Calabi flow and the corresponding hyperbolic combinatorial Calabi flow with surgery.
Please refer to Section \ref{section 4} for the definitions.
The main result for hyperbolic combinatorial Calabi flow with surgery is as follows.

\begin{theorem}\label{main theorem for hyper Calabi flow with surgery intro}
Suppose $(S, V)$ is a closed connected marked surface with $\chi(S)<0$
and $d_0$ is any initial piecewise hyperbolic metric on $(S, V)$.
Then the hyperbolic combinatorial Calabi flow with surgery
%\begin{equation}\label{hyperbolic CCF with surgery intro}
%\begin{aligned}
%\frac{du_i}{dt}=\Delta^\mathbb{H} F_i.
%\end{aligned}
%\end{equation}
exists for all time and
converges exponentially fast to a hyperbolic metric $d^*$ on $S$ after finite number of surgeries.
\end{theorem}

Theorem \ref{main theorem for hyper Calabi flow with surgery intro} is
an analogue of the results obtained in \cite{CXX, CP} with initial metric given by a piecewise hyperbolic metric,
which is a hyperbolic cone metric.
In Section \ref{section 4}, we will prove a generalization of Theorem \ref{main theorem for hyper Calabi flow with surgery intro}.

The main idea of the paper comes from reading of \cite{Ge-thesis, Ge, GGLSW, GLSW, L1}.
As the main tools used in this paper come from \cite{GGLSW, GLSW},
many notations are taken from \cite{GGLSW, GLSW} for consistence.

\subsection{Notations and definitions}
Here we give some notations and definitions used in the main results.
Suppose $S$ is a closed surface and $V$ is a finite subset of $S$, $(S, V)$ is called a marked surface.
A piecewise linear metric on $(S, V)$ is a flat cone metric with cone points contained in $V$.
The combinatorial curvature $K_i$ at $v_i\in V$ is $2\pi$ less the cone angle at $v_i$.
Suppose $\mathcal{T}=(V, E, F)$ is a triangulation of $(S, V)$, where $V,E,F$ represent the set of vertices, edges and faces respectively.
We use $(S, V, \mathcal{T})$ to denote a marked surface $(S, V)$ with a fixed triangulation $\mathcal{T}$.
If a map $l: E\rightarrow (0, +\infty)$ satisfies that $l_{rs}<l_{rt}+l_{st}$ for $\{r,s,t\}=\{i,j,k\}$, where $\{i,j,k\}$ is any triangle in $F$,
then $l$ determines a piecewise linear metric on $(S, V)$.
Given $(S, V)$ with a triangulation $\mathcal{T}$ and a map $l: E\rightarrow (0, +\infty)$ determined by a piecewise linear metric $d$ on $(S, V)$,
the vertex scaling \cite{L1,RW} of $d$
by a function $u: V\rightarrow \mathbb{R}$ is defined to be the piecewise linear metric $u*d$ on $(S, V)$
determined by $u*l: E\rightarrow (0, +\infty)$  with
\begin{equation*}
(u*l)_{ij}:=e^{u_i+u_j}l_{ij}, \ \ \forall \{ij\}\in E.
\end{equation*}
The function $u: V\rightarrow \mathbb{R}$ is called a conformal factor. Note that the vertex scaling of a piecewise linear metric on a surface
depends on the triangulation of the surface.

The Calabi flow (\ref{Calabi flow}) on smooth Riemannian surfaces deforms the metrics conformally.
If $g$ is conformal to a background metric $g_0$ with $g=e^ug_0$, where $u$ is the conformal factor,
then the Calabi flow (\ref{Calabi flow}) translates into an evolution equation
$$\frac{du}{dt}=\Delta_gK$$
of the conformal factor.
The Euclidean combinatorial Calabi flow for piecewise linear metrics is defined similarly.
\begin{definition}[\cite{Ge-thesis} Definition 2.11]\label{intro Euclidean Comb Calabi flow with trian}
Suppose $d_0$ is a piecewise linear metric on a triangulated surface $(S, V, \mathcal{T})$.
The Euclidean combinatorial Calabi flow on $(S, V, \mathcal{T})$ is defined as
\begin{equation}\label{Euclidean Calabi flow for triangulated surface intro}
\begin{aligned}
\left\{
  \begin{array}{ll}
    \frac{du_i}{dt}=\Delta^{\mathbb{E},\mathcal{T}} K_i & \hbox{ } \\
    u_i(0)=0 & \hbox{ }
  \end{array},
\right.
\end{aligned}
\end{equation}
where $u: V\rightarrow \mathbb{R}$ is the conformal factor and
$\Delta^{\mathbb{E},\mathcal{T}}$ is the Euclidean discrete Laplace operator of $u*d_0$
on $(S, V, \mathcal{T})$ defined as
\begin{equation}\label{laplace operator}
\begin{aligned}
(\Delta^{\mathbb{E}, \mathcal{T}} f)_i=\sum_{j; j\sim i} \omega_{ij}^{\mathbb{E}, \mathcal{T}}(f_j-f_i), \ \forall f\in \mathbb{R}^V
\end{aligned}
\end{equation}
with
$$\omega^{\mathbb{E}, \mathcal{T}}_{ij}=\cot \theta_{k}^{ij}+\cot\theta_l^{ij}.$$
Here $\theta_{i}^{jk}$ is the inner angle at the vertex $v_i$ in a triangle $\triangle ijk\in F$.
\end{definition}

The Euclidean discrete Laplace operator $\Delta^{\mathbb{E},\mathcal{T}}$ is the well-known finite-elements Laplace operator.
The Euclidean combinatorial Calabi flow (\ref{Euclidean Calabi flow for triangulated surface intro})
is defined for a fixed triangulation $\mathcal{T}$
of $(S,V)$ and may develop singularities,
including the conformal factor tends to infinity and some triangle degenerates along the flow.
To handle the possible singularities along the flow, we do surgery on the flow by flipping,
the idea of which comes from \cite{GGLSW,GLSW, L1}.
Note that the weight $\omega^{\mathbb{E}, \mathcal{T}}_{ij}$ may be negative and
$\omega^{\mathbb{E}, \mathcal{T}}_{ij}\geq 0$ if and only if $\theta_{k}^{ij}+\theta_l^{ij}\leq \pi$, which is the locally Delaunay condition
of $\mathcal{T}$ on the edge $\{ij\}$ \cite{BS}.
To ensure that the Euclidean discrete Laplace operator $\Delta^{\mathbb{E},\mathcal{T}}$ has good properties,
we require the triangulations along the
Euclidean combinatorial Calabi flow (\ref{Euclidean Calabi flow for triangulated surface intro}) to be Delaunay, which is equivalent to
every edge satisfies the locally Delaunay condition \cite{BS}.
Note that each piecewise linear metric on $(S, V)$ has at least one Delaunay triangulation \cite{AK, BS},
so this additional condition is reasonable.
Along the Euclidean combinatorial Calabi flow (\ref{Euclidean Calabi flow for triangulated surface intro})
on $(S, V)$ with a triangulation $\mathcal{T}$, if $\mathcal{T}$ is Delaunay in $u(t)*d_0$  for
$t\in [0, T]$ and not Delaunay in $u(t)*d_0$  for $t\in (T, T+\epsilon)$, $\epsilon>0$, there exists
an edge $\{ij\}\in E$ such that $\theta_{k}^{ij}(t)+\theta_l^{ij}(t)\leq \pi$
for $t\in [0, T]$ and $\theta_{k}^{ij}(t)+\theta_l^{ij}(t)> \pi$ for $t\in (T, T+\epsilon)$.
We replace the triangulation $\mathcal{T}$ by a new triangulation $\mathcal{T}'$ at time $t=T$
by replacing two triangles $\triangle ijk$ and $\triangle ijl$
adjacent to $\{ij\}$ by two new triangles $\triangle ikl$ and $\triangle jkl$.
This is called a \textbf{surgery by flipping} on the triangulation $\mathcal{T}$, which is an isometry of $(S, V)$ in the piecewise linear metric $u(T)*d_0$.
After the surgery at time $t=T$, we run the Euclidean combinatorial Calabi flow (\ref{Euclidean Calabi flow for triangulated surface intro})
on $(S, V, \mathcal{T}')$ with initial metric coming from
the Euclidean combinatorial Calabi flow (\ref{Euclidean Calabi flow for triangulated surface intro}) on $(S, V, \mathcal{T})$ at time $t=T$.

Vertex scaling, combinatorial Calabi flow and surgery by flipping can also defined for piecewise hyperbolic metrics on surfaces.
Please refer to Section \ref{section 2} and Section \ref{section 4} for more details.
%\begin{definition}\label{defn of Euclid CCF with surgery}
%The combinatorial Calabi flow with surgery is
%\begin{equation}\label{Euclidean CCF with surgery intro}
%\begin{aligned}
%\frac{du_i}{dt}=\Delta^\mathbb{E} F_i,
%\end{aligned}
%\end{equation}
%where $\Delta^\mathbb{E}$ is the Laplace operator with respect to some Delaunay triangulation in PL metric $u(t)*d$.
%\end{definition}
%The main result of the paper is the following result on Euclidean combinatorial Calabi flow with surgery.
%\begin{theorem}\label{main theorem for Calabi flow with surgery intro}
%Suppose $(S, V)$ is a closed connected marked surface and $d_0$ is any initial piecewise linear metric on $(S, V)$.
%Then the combinatorial Calabi flow with surgery
%exists for all time and
%converges exponentially fast to a piecewise linear metric $d^*$ with
%constant combinatorial curvature after finite number of surgeries.
%\end{theorem}

\subsection{Organization of the paper}
The paper is organized as follows.
In Section \ref{section 2}, we give some preliminaries on discrete conformal geometry,
including the definitions of polyhedral metrics, vertex scaling, discrete curvature and discrete Laplace operators.
In Section \ref{section 3}, we study the Euclidean combinatorial Calabi flow on surfaces
and prove a generalization of Theorem \ref{main theorem for Calabi flow with surgery intro}.
In Section \ref{section 4}, we study the hyperbolic combinatorial Calabi flow on surfaces
and prove a generalization of Theorem \ref{main theorem for hyper Calabi flow with surgery intro}.
In Section \ref{section 5}, we give some remarks and propose an interesting question.

\section{Polyhedral metrics, discrete curvatures and discrete Laplace operators}\label{section 2}
\subsection{Polyhedral metrics on surfaces}
The definition of piecewise linear metric on marked surfaces has been given in Section \ref{section 1}.
Here we extend the definition to polyhedral metrics on surfaces.
\begin{definition}[\cite{GGLSW, GLSW}]
Suppose $(S, V)$ is a marked surface.
A piecewise linear (hyperbolic and spherical respectively) metric on  $(S, V)$
is a flat (hyperbolic and spherical respectively) cone metric on $S$ whose
cone points are contained in $V$.
Piecewise linear, piecewise hyperbolic and piecewise spherical metrics on marked surfaces are all called polyhedral metrics.
\end{definition}

In this paper, we concern only about piecewise linear and piecewise hyperbolic metrics.
For simplicity, piecewise linear and piecewise hyperbolic metrics
are denoted by PL and PH metrics respectively in the following.

A marked surface $(S, V)$ with a triangulation $\mathcal{T}=\{V,E,F\}$ is called a triangulated surface
and denoted by $(S, V, \mathcal{T})$.
In this paper, a function defined on vertices $V$ is regarded as a column vector and $n=\#V$ is used to denote the number of vertices. Moreover, all vertices, marked by $v_{1},...,v_{n}$, are supposed to be ordered one by one and we often write $i$ instead of $v_i$.
We use $\{ij\}$ to denote the edge between $v_i$ and $v_j$ in E and use $\triangle ijk$ to denote the face determined by $v_i, v_j$ and $v_k$ in $F$.

Geometrically, PL metrics on $(S, V)$ are obtained by isometrically gluing Euclidean triangles along their edges so that
the cone points are contained in $V$, which gives a geometric triangulation $\mathcal{T}$ of $(S, V)$ whose
simplices are quotients of the simplices in the disjoint union.
A triangulation $\mathcal{T}$ of $(S, V)$ is geometric in a PL
metric $d$ on $(S, V)$ if each triangle in $\mathcal{T}$ is isometric to a Euclidean triangle in $d$.
PH metrics on $(S, V)$ are obtained similarly with Euclidean triangles replaced by hyperbolic triangles.
A geometrical triangulation of $(S, V)$ with a PL metric is said to be a \emph{Delaunay triangulation}
if the sum of two angles facing each edge is at most $\pi$.
A geometrical triangulation of $(S, V)$ with a PH metric is said to be a \emph{Delaunay triangulation}
if for each edge $e$ adjacent to two hyperbolic triangles $t$ and $t'$ , the interior of the circumball
of $t$ does not contain the vertices of $t'$ when the quadrilateral $t\cup t'$ is lifted to $\mathbb{H}^2$.
See \cite{AK, BS, GGLSW, GLSW, Lei, R} for further discussions on Delaunay triangulations of surfaces.
Note that a polyhedral metric on a marked surface is independent of the triangulations.

Suppose $\mathcal{T}=\{V,E,F\}$ is a geometric triangulation of $(S, V)$ with a PL or PH metric $d$,
then the metric $d$ determines a map
\begin{equation*}
\begin{aligned}
d: E&\longrightarrow (0, +\infty)\\
   \{ij\}&\mapsto d_{ij}\triangleq d(\{ij\})
\end{aligned}
\end{equation*}
such that for any topological triangle $\triangle ijk\in F$, the triangle inequalities
$d_{ij}<d_{ik}+d_{jk}, d_{ik}<d_{ij}+d_{jk}, d_{jk}<d_{ij}+d_{ik}$
are satisfied. Conversely, given a map $d: E\rightarrow (0, +\infty)$ satisfying the triangle inequalities
$d_{ij}<d_{ik}+d_{jk}, d_{ik}<d_{ij}+d_{jk}, d_{jk}<d_{ij}+d_{ik}$
for each triangle $\triangle ijk\in F$, the map $d: E\rightarrow (0, +\infty)$ uniquely determines a PL metric on $(S, V)$
by the construction of the polyhedral metrics.
Then for a triangulated surface $(S, V, \mathcal{T})$, the space of PL metrics is given by
\begin{equation*}
\begin{aligned}
\mathbb{R}^{E(\mathcal{T})}_{\triangle}=\{d\in \mathbb{R}^{E(\mathcal{T})}_{>0}| d_{ij}, d_{ik}, d_{jk}
&\text{ satisfy the triangle}\\
&\text{ inequalities for any } \triangle ijk\in F\}.
\end{aligned}
\end{equation*}
And there is an injective map
$\Phi_{\mathcal{T}}: \mathbb{R}^{E(\mathcal{T})}_{\triangle}\longrightarrow T_{PL}(S, V)$
between the subspace
$\mathbb{R}^{E(\mathcal{T})}_{\triangle}$
and the Teichim\"{u}ller space $T_{PL}(S, V)$ of PL metrics on $(S, V)$.
Note that $\mathbb{R}^{E(\mathcal{T})}_{\triangle}$ is a proper subset of $\mathbb{R}^{E(\mathcal{T})}_{>0}$.

\subsection{Discrete curvature}

On a marked surface, the well-known combinatorial curvature is defined as follows.

\begin{definition}
Suppose $(S, V)$ is a marked surface with a polyhedral metric,
the combinatorial curvature $K_i$ at $v_i\in V$ is $2\pi$ less the cone angle at $v_i$.
\end{definition}

If $\mathcal{T}$ is a geometric triangulation of $(S, V)$,
the combinatorial curvature $K_i$ at  $v_i$ is
\begin{equation*}
K_i=2\pi-\sum_{\triangle ijk\in F} \theta_i^{jk},
\end{equation*}
where the summation is taken over triangles with $v_i$ as vertex and $\theta_i^{jk}$ is the inner angle
of $\triangle ijk$ at $v_i$.
Note that the combinatorial curvature $K_i$ is independent of the geometric triangulations of $(S, V)$ with a given polyhedral metric.
Combinatorial curvature $K$ satisfies the following discrete Gauss-Bonnet formula \cite{CL}
\begin{equation*}
\sum_{i=1}^n K_i=2\pi\chi(S)-\lambda Area(S),
\end{equation*}
where $Area(S)$ is the area of the marked surface $(S, V)$ with polyhedral metric $d$ and $\lambda=-1, 0, +1$ respectively
when $d$ is a hyperbolic, Euclidean or spherical polyhedral metric respectively.

\subsection{Vertex scaling of polyhedral metrics}
Vertex scaling of PL metrics on a triangulated surface was introduced
by Luo \cite{L1} and R\v{o}cek-Williams \cite{RW} independently as an analogy of
the conformal transformation of Riemannian metrics.

\begin{definition}[\cite{L1, RW}]
Suppose $d$ is a PL metric on a triangulated surface $(S, V, \mathcal{T})$ and
$u$ is a function defined on the vertices $V$. The vertex scaling of $d$ by $u$ for $(S, V, \mathcal{T})$
is defined to be the PL metric $u*d$ such that
\begin{equation*}
(u*d)_{ij}:=d_{ij}e^{u_i+u_j}
\end{equation*}
determines a PL metric $u*d$ on $(S, V, \mathcal{T})$. $u$ is called a conformal factor.
\end{definition}
Note that $u*d$ determines a PL metric on a triangulated surface $(S, V, \mathcal{T})$ is equivalent to $u: V\rightarrow \mathbb{R}$
is in the following admissible space of the conformal factors
\begin{equation*}
\Omega^{\mathbb{E}, \mathcal{T}}(d)\triangleq \{u\in \mathbb{R}^V| d_{rs}e^{u_t}+d_{rt}e^{u_s}>d_{st}e^{u_r}, \{r,s,t\}=\{i,j,k\}, \forall\triangle ijk\in F \}.
\end{equation*}

Luo-Sun-Wu \cite{LSW} showed that the vertex scaling of
the PL metric is an approximation of the conformal transformation in the smooth case.

Vertex scaling of PH metrics on a triangulated surface was introduced by Bobenko-Pinkall-Springborn \cite{BPS}.
\begin{definition}[\cite{BPS}]
Suppose $d$ is a PH metric on a triangulated surface $(S, V, \mathcal{T})$ and
$u$ is a function defined on the vertices $V$. The vertex scaling of $d$ by $u$
is defined to be $u*d$ such that
\begin{equation*}
\sinh\frac{(u*d)_{ij}}{2}:=e^{u_i+u_j}\sinh\frac{d_{ij}}{2}
\end{equation*}
determines a PH metric $u*d$ on $(S, V, \mathcal{T})$. $u$ is called a conformal factor.
\end{definition}

Similar to the PL metrics, $u*d$ determines a PH metric
if and only if the conformal factor $u$ is in the admissible space $\Omega^{\mathbb{H}, \mathcal{T}}(d)$ of conformal factors
\begin{equation*}
\begin{aligned}
\Omega^{\mathbb{H}, \mathcal{T}}(d)\triangleq \{u\in \mathbb{R}^V| (u*d)_{rs}+&(u*d)_{rt}>(u*d)_{st},\\
&\{r,s,t\}=\{i,j,k\}, \forall\triangle ijk\in F \}.
\end{aligned}
\end{equation*}

The combinatorial curvature $K$ is rigid with respect to the conformal factor.
\begin{theorem}[\cite{BPS}]\label{global rigidity for fix triangulation}
Suppose $(S, V, \mathcal{T})$ is a triangulated surface with a PL or PH metric, then
the conformal factor is uniquely determined by the combinatorial curvature $K$ (up to scaling for the PL metric).
\end{theorem}

Vertex scaling for spherical polyhedral metrics was defined by Bobenko-Pinkall-Springborn \cite{BPS}. Gu-Luo-Sun-Wu \cite{GLSW} observed
that the notions of vertex scaling are related to the Ptolemy identities for all polyhedral surfaces.
Note that the definition of vertex scaling for polyhedral metrics depends on
the triangulation of the marked surface $(S, V)$.

\subsection{Laplace operators on triangulated surfaces}

The discrete Laplace operator of a PL metric on a triangulated surface, known as finite elements Laplacian,
has been extensively studied in geometry and computer graphics and is defined as follows.

\begin{definition}\label{Euclidean Laplacian with triang}
Suppose $(S, V, \mathcal{T})$ is a triangulated surface with a PL metric $d$.
The Euclidean discrete Laplace operator of $d$ on $(S, V, \mathcal{T})$ is defined to be the map
\begin{equation*}
\begin{aligned}
\Delta^{\mathbb{E}, \mathcal{T}}: \mathbb{R}^V&\longrightarrow \mathbb{R}^V\\
f&\mapsto \Delta^{\mathbb{E}, \mathcal{T}} f,
\end{aligned}
\end{equation*}
where $f: V\rightarrow \mathbb{R}$ is a function defined on the vertices
and the value of $\Delta^{\mathbb{E}, \mathcal{T}} f$ at $v_i$ is
\begin{equation}\label{laplace operator}
\begin{aligned}
(\Delta^{\mathbb{E}, \mathcal{T}} f)_i=\sum_{j; j\sim i} \omega_{ij}^{\mathbb{E}, \mathcal{T}}(f_j-f_i)
\end{aligned}
\end{equation}
with weight
$$\omega^{\mathbb{E}, \mathcal{T}}_{ij}=\cot \theta_{k}^{ij}+\cot\theta_l^{ij}.$$
Here $\theta_{i}^{jk}$ is the inner angle at the vertex $v_i$ in a triangle $\triangle ijk\in F$.
\end{definition}

\begin{remark}\label{intrinsic of Euclid Laplace}
Generally, for a triangulated surface $(S, V, \mathcal{T})$ with a PL metric $d$,
the Euclidean discrete Lapalce operator of $d$ depends on the triangulation $\mathcal{T}$.
However, if the triangulation $\mathcal{T}$ is Delaunay, even though there exists different
Delaunay triangulations for the same PL metric $d$ on a marked surface $(S, V)$,
the Euclidean discrete Laplace operator $\Delta^{\mathbb{E}, \mathcal{T}}$
is independent of the Delaunay triangulations of $d$ on $(S, V)$ \cite{BS}.
In this case, the Euclidean discrete Laplace operator $\Delta^{\mathbb{E}, \mathcal{T}}$ is intrinsic
in the sense that it is
independent of the Delaunay triangulations.
\end{remark}

Suppose $u: V\rightarrow \mathbb{R}$ is a conformal factor defined on the vertices and
the vertex $v_j$ is adjacent to $v_i$, Luo \cite{L1} proved
\begin{equation*}
\frac{\partial \theta_{i}^{jk}}{\partial u_j}=\cot\theta_k^{ij},\ \
\frac{\partial K_i}{\partial u_j}=-(\cot\theta_k^{ij}+\cot\theta_l^{ij}),
\end{equation*}
where $\triangle ijk$ and $\triangle ijl$ are adjacent triangles in $F$.

Set
\begin{equation}\label{definition of L}
\begin{aligned}
L=(L_{ij})_{N\times N}=\frac{\partial (K_1, \cdots, K_N)}{\partial (u_1, \cdots, u_N)}
=\left(
                                                                                         \begin{array}{ccc}
                                                                                           \frac{\partial K_1}{\partial u_1} & \cdots & \frac{\partial K_1}{\partial u_N} \\
                                                                                           \vdots & \ddots & \vdots \\
                                                                                           \frac{\partial K_N}{\partial u_1} & \cdots & \frac{\partial K_N}{\partial u_N} \\
                                                                                         \end{array}
                                                                                       \right)
\end{aligned}
\end{equation}
be the Jacobian of $K$ with respect to $u$.
For Euclidean polyhedral metrics, the matrix $L$ has the following property.

\begin{lemma}[\cite{L1}]\label{property of Euclidean L}
For a triangulated surface $(S, V, \mathcal{T})$ with a PL metric $d$,
the matrix $L$ is symmetric and positive semi-definite on $\Omega^{\mathbb{E}, \mathcal{T}}(d)$ with kernel $\{t \mathbf{1}|t\in \mathbb{R}\}$,
where $\mathbf{1}=(1, \cdots, 1)^T$.
\end{lemma}

By Lemma \ref{property of Euclidean L}, we have $\sum_{j=1}^N\frac{\partial K_i}{\partial u_j}=0$.
Then the Euclidean discrete Laplace operator (\ref{laplace operator}) translates into
\begin{equation*}
\begin{aligned}
(\Delta^{\mathbb{E}, \mathcal{T}} f)_i=\sum_{j; j\sim i}(-\frac{\partial K_i}{\partial u_j})(f_j-f_i)=-\sum_{j\in V}\frac{\partial K_i}{\partial u_j}f_j,
\end{aligned}
\end{equation*}
which implies that
$$\Delta^{\mathbb{E}, \mathcal{T}}=-L$$
in matrix form.
In this way, we can take the Euclidean discrete Laplace operator $\Delta^{\mathbb{E}, \mathcal{T}}$
as a matrix-valued map defined on the admissible space of conformal factors $\Omega^{\mathbb{E}, \mathcal{T}}(d)$.
If the triangulation $\mathcal{T}$ is Delaunay in $u*d$, then the
Euclidean discrete Laplace operator $\Delta^{\mathbb{E}, \mathcal{T}}$, as a matrix-valued map,
is defined on
\begin{equation*}
\Omega^{\mathbb{E}, \mathcal{T}}_{D}(d)\triangleq \{u\in \mathbb{R}^V| \mathcal{T} \text{ is Delaunay in } u*d  \},
\end{equation*}
which is a subspace of the admissible conformal factors $\Omega^{\mathbb{E}, \mathcal{T}}(d)$ \cite{GLSW,Penner}.
The idea of writing the discrete Laplace operator in matrix form comes from Ge \cite{Ge-thesis, Ge}.

Note that $L$ is the Jacobian of the curvature $K$ with respect to the conformal factor $u$,
we define the hyperbolic discrete Laplace operator similarly.
\begin{definition}
Suppose $(S, V, \mathcal{T})$ is a triangulated surface with a PH metric $d$.
The hyperbolic discrete Laplace operator is defined to be the map
\begin{equation*}
\begin{aligned}
\Delta^{\mathbb{H}, \mathcal{T}}: \mathbb{R}^V&\longrightarrow \mathbb{R}^V\\
f&\mapsto \Delta^{\mathbb{H}, \mathcal{T}} f,
\end{aligned}
\end{equation*}
where $f: V\rightarrow \mathbb{R}$ is a function defined on the vertices
and the value of $\Delta^{\mathbb{H}, \mathcal{T}} f$ at $v_i$ is
\begin{equation}\label{hyper laplace operator}
\begin{aligned}
(\Delta^{\mathbb{H}, \mathcal{T}} f)_i=-(Lf)_i
\end{aligned}
\end{equation}
with $L$ formally given by $(\ref{definition of L})$.
\end{definition}

For hyperbolic polyhedral metrics, the matrix $L$ has the following property.

\begin{lemma}[\cite{BPS}, Proposition 6.1.5]
For a triangulated surface $(S, V, \mathcal{T})$ with a PH metric $d$,
the matrix $L$ is symmetric and strictly positive definite on $\Omega^{\mathbb{H}, \mathcal{T}}(d)$.
\end{lemma}

Similar to the Euclidean discrete Laplace operator $\Delta^{\mathbb{E}, \mathcal{T}}$,
the hyperbolic discrete Laplace operator $\Delta^{\mathbb{H}, \mathcal{T}}$
is a matrix-valued map defined on the admissible space of conformal factors $\Omega^{\mathbb{H}, \mathcal{T}}(d)$.
If the triangulation $\mathcal{T}$ is Delaunay in $u*d$, then the
hyperbolic discrete Laplace operator $\Delta^{\mathbb{E}, \mathcal{T}}$, as a matrix-valued map,
is defined on
\begin{equation*}
\Omega^{\mathbb{H}, \mathcal{T}}_{D}(d)\triangleq \{u\in \mathbb{R}^V| \mathcal{T} \text{ is Delaunay in } u*d  \},
\end{equation*}
which is a subspace of the admissible conformal factors $\Omega^{\mathbb{H}, \mathcal{T}}(d)$ \cite{GGLSW}.

\begin{remark}
Generally, for a triangulated surface $(S, V, \mathcal{T})$ with a PH metric $d$, the hyperbolic discrete Lapalce operator
$\Delta^{\mathbb{H}, \mathcal{T}}$
depends on the triangulations.
However, if the triangulation $\mathcal{T}$ is restricted to be Delaunay,
the hyperbolic discrete Laplace operator $\Delta^{\mathbb{H}, \mathcal{T}}$ is independent of
the Delaunay triangulations of $(S, V)$ for the PH metric $d$ \cite{GGLSW}.
In this case, the hyperbolic discrete Laplace operator by $\Delta^{\mathbb{H}, \mathcal{T}}$
is intrinsic in the sense that it is
independent of the Delaunay triangulations.
\end{remark}

\section{Euclidean combinatorial Calabi flow}\label{section 3}
\subsection{Euclidean combinatorial Calabi flow on triangulated surfaces}
The definition of Euclidean combinatorial Calabi flow on a triangulated marked surface $(S, V, \mathcal{T})$
is given in Definition \ref{intro Euclidean Comb Calabi flow with trian}, which could be written in the following
matrix form \cite{Ge-thesis}
\begin{equation*}
\frac{du}{dt}=-LK.
\end{equation*}
Note that this is essential an ODE system, therefor the Euclidean combinatorial Calabi flow
exists in a short time \cite{Ge-thesis}.
%And when the Euclidean combinatorial Calabi flow exists, the conformal factors lie in
%the admissible space of conformal factors $\Omega^{\mathbb{E}, \mathcal{T}}(d_0)$.
By direct calculations, we have the combinatorial curvature $K$ evolves
according to
\begin{equation*}
\frac{dK}{dt}=-(\Delta^{\mathbb{E},\mathcal{T}})^2 K
\end{equation*}
along the Euclidean combinatorial Calabi flow (\ref{Euclidean Calabi flow for triangulated surface intro}).
By Proposition \ref{property of Euclidean L},
$\sum_{i=1}^nu_i$ is invariant along the Euclidean combinatorial Calabi flow (\ref{Euclidean Calabi flow for triangulated surface intro}) \cite{Ge-thesis}.

Similar to the fact that combinatorial Calabi flow for Thurston's circle packing metrics is
the negative gradient flow of combinatorial Calabi energy \cite{Ge-thesis, Ge, GX2},
the combinatorial Calabi flow (\ref{Euclidean Calabi flow for triangulated surface intro}) for
PL metrics is also the negative gradient flow of combinatorial Calabi energy.

\begin{definition}[\cite{Ge-thesis, Ge}]
Suppose $d$ is a polyhedral metric on a marked surface $(S, V)$, the combinatorial Calabi energy of $d$ on $(S, V)$
is defined to be
\begin{equation*}
\mathcal{C}=||K||^2=\sum_{i=1}^nK_i^2.
\end{equation*}
If $K^*$ is a function defined on $V$, the modified combinatorial Calabi energy is defined to be
\begin{equation*}
\overline{\mathcal{C}}=||K-K^*||^2=\sum_{i=1}^n(K_i-K_i^*)^2.
\end{equation*}
\end{definition}

\begin{remark}
Combinatorial Calabi energy $\mathcal{C}$ and modified combinatorial Calabi energy $\overline{\mathcal{C}}$ are independent of
the geometric triangulations of $d$ on $(S, V)$.
\end{remark}

Note that
$$\nabla_u \mathcal{C}=2LK=-2\Delta^{\mathbb{E}, \mathcal{T}} K.$$
The Euclidean combinatorial Calabi flow (\ref{Euclidean Calabi flow for triangulated surface intro})  can be written as
$$\frac{du}{dt}=-\frac{1}{2}\nabla_u \mathcal{C}.$$
\begin{proposition}[\cite{Ge-thesis}]
The Euclidean combinatorial Calabi flow $(\ref{Euclidean Calabi flow for triangulated surface intro})$
on a triangulated surface $(S, V, \mathcal{T})$
is the negative gradient flow of combinatorial Calabi energy $\mathcal{C}$
and the combinatorial Calabi energy $\mathcal{C}$ is decreasing along the flow.
\end{proposition}

Similar to the stability results in \cite{GX1,GX2,GX3,GX4,GX5}, we have the following result for
Euclidean combinatorial Calabi flow (\ref{Euclidean Calabi flow for triangulated surface intro}).

\begin{theorem}\label{stability for triangulated CCF}
If the solution of Euclidean combinatorial Calabi flow $(\ref{Euclidean Calabi flow for triangulated surface intro})$
on a triangulated surface $(S, V, \mathcal{T})$
converges as time tends to infinity, then the limit metric is a constant combinatorial curvature PL metric.
Furthermore, suppose there is a constant combinatorial curvature PL metric $d^*=u^**d_0$
on $(S, V, \mathcal{T})$ with $\sum_{i=1}^nu^*_{i}=0$,
there exists a constant $\delta>0$
such that if the initial modified combinatorial Calabi energy $||K(u(0))-K(u^*)||^2<\delta$,
then the Euclidean combinatorial Calabi flow $(\ref{Euclidean Calabi flow for triangulated surface intro})$
on $(S, V, \mathcal{T})$ exists for all time and converges exponentially fast to $u^*$.
\end{theorem}

\proof
If $u(t)$ converges as time tends to infinity,
then $u(+\infty)=\lim_{t\rightarrow +\infty} u(t)\in \Omega^{\mathbb{E}, \mathcal{T}}(d_0)$ exists.
As $K$ is a smooth function of $u\in \Omega^{\mathbb{E}, \mathcal{T}}(d_0)$, we have
$K(+\infty)=\lim_{t\rightarrow +\infty}K(u(t))$ and $L(+\infty)=\lim_{t\rightarrow +\infty}L(u(t))$ exist.
Similarly, $\mathcal{C}(+\infty)$ and $\mathcal{C}'(+\infty)$ exist.
Note that
$$\mathcal{C}'(t)=-2K^TL^2K\leq 0$$
and $\mathcal{C}(t)$ is uniformly bounded, we have
$$\mathcal{C}'(+\infty)=-2K^T(+\infty)L^2(+\infty)K(+\infty)=0.$$
By Lemma \ref{property of Euclidean L}, we have $K(+\infty)=\frac{2\pi\chi(S)}{n}(1, 1, \cdots, 1)^T$,
which implies that $u(+\infty)*d_0$ is a PL metric with constant combinatorial curvature.

By Theorem \ref{global rigidity for fix triangulation}, we know that $u^*=u(+\infty)$ is the unique conformal factor
such that $u^**d_0$ has constant combinatorial curvature.
Set $\Gamma(u)=\Delta^{\mathbb{E},\mathcal{T}} K=-LK$, then $\Gamma(u^*)=0$ and
\begin{equation}\label{negative definiteness of critical point}
\begin{aligned}
D_u\Gamma(u^*)=-L^TL\leq 0.
\end{aligned}
\end{equation}
Note that
$rank\ (D_{u}\Gamma(u^*))=n-1$ and the kernel of $D_{u}\Gamma(u^*)$ is
exactly $c(1, 1, \cdots, 1)^T,\,c\in \mathbb{R}$.
Along the Euclidean combinatorial Calabi flow (\ref{Euclidean Calabi flow for triangulated surface intro}),
$\sum_{i=1}^nu_i$ is a constant.
This implies that $u^*$ is a local attractor of the Euclidean combinatorial Calabi flow (\ref{Euclidean Calabi flow for triangulated surface intro})
and the conclusion follows from the Lyapunov Stability Theorem (\cite{P}, Chapter 5). \qed

Theorem \ref{stability for triangulated CCF} does not give the long time existence and convergence of the Euclidean
combinatorial Calabi flow (\ref{Euclidean Calabi flow for triangulated surface intro}) for general initial PL metrics.
To study the long time behavior of the flow for general initial PL metrics,
we need to analysis the Euclidean discrete Laplace operator $\Delta^{\mathbb{E}, \mathcal{T}}$.
Note that, for $j\sim i$, the weight of the Euclidean discrete Laplace operator $\Delta^{\mathbb{E}, \mathcal{T}}$ is
\begin{equation*}
\begin{aligned}
\omega^{\mathbb{E}, \mathcal{T}}_{ij}
=-\frac{\partial K_i}{\partial u_j}=\cot \theta_{k}^{ij}+\cot\theta_l^{ij}
=\frac{\sin(\theta_{k}^{ij}+\theta_{l}^{ij})}{\sin\theta_{k}^{ij}\sin\theta_{l}^{ij}},
\end{aligned}
\end{equation*}
which may be negative and unbounded along the Euclidean combinatorial
Calabi flow (\ref{Euclidean Calabi flow for triangulated surface intro}) for general initial PL metric.
Even through discrete Laplace operators with negative coefficients has many applications in
the study of combinatorial curvatures on two and three dimensional manifolds (see \cite{CR, GX5, G1,G2,Guo,L2,R1,X1,X2,Z} for example),
a discrete Laplace operator is generally defined with nonnegative weights \cite{CFR}.
For $\Delta^{\mathbb{E}, \mathcal{T}}$, the weight is nonnegative is equivalent to
$\theta_{k}^{ij}+\theta_{l}^{ij}\leq \pi,$
which is the Delaunay condition \cite{BS}.
So it is natural to require the triangulation to be Delaunay, otherwise
the coefficients of the discrete Laplace operator will be negative.
Note that the admissible space for a Delaunay triangulation
\begin{equation*}
\begin{aligned}
D_{PL}(\mathcal{T})=\{[d]\in T_{PL}(S, V)| \mathcal{T} \text{ is isotopic to a Delaunay triangulation of } d\}
\end{aligned}
\end{equation*}
is a subspace of the admissible space
\begin{equation*}
\begin{aligned}
\Omega(\mathcal{T})=\{[d]| d_{ij}, d_{ik}, d_{jk} \text{ satisfy the triangle inequalities for any } \triangle ijk\in F\}.
\end{aligned}
\end{equation*}

\subsection{Gu-Luo-Sun-Wu's work on discrete uniformization theorem}
One of the main tools used in the proof of Theorem \ref{main theorem for Calabi flow with surgery intro} for convergence of combinatorial Calabi flow
with surgery is
the discrete conformal theory developed by Gu-Luo-Sun-Wu \cite{GLSW}.
In this subsection, we will briefly recall the theory. For details of the theory,
please refer to \cite{GLSW}.

\begin{definition}[\cite{GLSW} Definition 1.1]\label{GLSW's definition of edcf}
Two piecewise linear metrics $d, d'$ on $(S, V)$ are discrete conformal if
there exist sequences of PL metrics $d_1=d$, $\cdots$,  $d_m=d'$
on $(S, V)$ and triangulations $\mathcal{T}_1, \cdots, \mathcal{T}_m$ of
$(S, V)$ satisfying
\begin{description}
  \item[(a)] (Delaunay condition) each $\mathcal{T}_i$ is Delaunay in $d_i$,
  \item[(b)] (Vertex scaling condition) if $\mathcal{T}_i=\mathcal{T}_{i+1}$, there exists a function
  $u:V\rightarrow \mathbb{R}$ so that if $e$ is an edge in $\mathcal{T}_i$ with end points $v$ and $v'$,
  then the lengths $l_{d_{i+1}}(e)$ and  $l_{d_{i}}(e)$ of $e$ in $d_i$ and $d_{i+1}$ are related by
  $$l_{d_{i+1}}(e)=l_{d_{i}}(e)e^{u(v)+u(v')},$$
  \item[(c)] if $\mathcal{T}_i\neq \mathcal{T}_{i+1}$, then $(S, d_i)$ is isometric to $(S, d_{i+1})$
  by an isometry homotopic to identity in $(S, V)$.
\end{description}
\end{definition}
The space of PL metrics on $(S, V)$ discrete conformal to $d$ is called the conformal class of $d$ and
denoted by $\mathcal{D}(d)$.

Recall the following result for Delaunay triangulations.
\begin{lemma}[\cite{AK, BS}]\label{Euclidean surgery}
If $\mathcal{T}$ and $\mathcal{T}'$ are Delaunay triangulations of $d$, then there exists
a sequence of Delaunay triangulations $\mathcal{T}_1=\mathcal{T}, \mathcal{T}_2, \cdots, \mathcal{T}_k=\mathcal{T}'$
so that $\mathcal{T}_{i+1}$ is obtained from $\mathcal{T}_i$ by a diagonal switch.
\end{lemma}

The diagonal switch in Lemma \ref{Euclidean surgery} is the surgery by flipping described in the introduction.

The following discrete uniformization theorem was established in \cite{GLSW}.

\begin{theorem}[\cite{GLSW} Theorem 1.2]\label{Euclidean discrete uniformization}
Suppose $(S, V)$ is a closed connected marked surface and $d$ is a PL metric on $(S, V)$.
Then for any $K^*: V\rightarrow (-\infty, 2\pi)$ with $\sum_{v\in V}K^*(v)=2\pi\chi(S)$,
there exists a PL metric $d'$, unique up to scaling and isometry homotopic to the identity
 on $(S, V)$, such that $d'$ is discrete conformal to $d$ and the discrete curvature
 of $d'$ is $K^*$.
Furthermore, the metric $d'$ can be found using a finite dimensional (convex) variational principle.
\end{theorem}

Denote the Teichim\"{u}ller space of all PL metrics on $(S, V)$ by $T_{PL}(S, V)$ and decorated
Teichim\"{u}ller space of all equivalence class of decorated hyperbolic metrics on $S-V$ by $T_D(S-V)$.
In the proof of Theorem \ref{Euclidean discrete uniformization}, Gu-Luo-Sun-Wu proved the following result.

\begin{theorem}[\cite{GLSW}]\label{main result of GLSW}
There is a $C^1$-diffeomorphism $\mathbf{A}: T_{PL}(S, V)\rightarrow T_D(S, V)$ between $T_{PL}(S, V)$ and $T_D(S-V)$.
Furthermore, the space $\mathcal{D}(d)\subset T_{PL}(S, V)$ of all equivalence classes of PL metrics
discrete conformal to $d$ is $C^1$-diffeomorphic to $\{p\}\times \mathbb{R}^V_{>0}$ under the diffeomorphism $\mathbf{A}$,
where $p$ is the unique hyperbolic metric on $S-V$ determined by the PL metric $d$ on $(S, V)$.
\end{theorem}

Set $u_i=\ln w_i$ for $w=(w_1, w_2, \cdots, w_n)\in \mathbb{R}^n_{>0}$.
Using the map $\mathbf{A}$, Gu-Luo-Sun-Wu \cite{GLSW} defined the curvature map
\begin{equation*}
\begin{aligned}
\mathbf{F}:\mathbb{R}^n&\rightarrow (-\infty, 2\pi)^n\\
u&\mapsto K_{\mathbf{A}^{-1}(p, w(u))}
\end{aligned}
\end{equation*}
and proved the following property of $\mathbf{F}$.

\begin{proposition}[\cite{GLSW}]\label{energy function Euclidean}
\begin{enumerate}
  \item For any $k\in \mathbb{\mathbb{R}}$, $\mathbf{F}(v+k(1, 1, \cdots, 1))=\mathbf{F}(v)$.
  \item There exists a $C^2$-smooth convex function $W=\int \sum_{i=1}^n \mathbf{F}_i(u)du_i: \mathbb{R}^n\rightarrow \mathbb{R}$
so that its gradient $\nabla W$ is $\mathbf{F}$ and the restriction $W: \{u\in \mathbb{R}^n|\sum_{i=1}^nu_i=0\}\rightarrow \mathbb{R}$
is strictly convex.
\end{enumerate}
\end{proposition}

%The convex function $W$ is in fact given by $W(u)=\int \sum_{i=1}^n \mathbf{F}_i(u)du_i$.

Theorem \ref{main result of GLSW} implies that the union of the admissible spaces $\Omega^{\mathbb{E}, \mathcal{T}}_{D}(d')$
of conformal factors such that $\mathcal{T}$ is Delaunay for $d'\in \mathcal{D}(d)$
is $\mathbb{R}^n$.
It is further proved that $\mathbb{R}^n=\cup_{\mathcal{T}} \Omega^{\mathbb{E}, \mathcal{T}}_{D}(d')$ is an analytic cell decomposition of $\mathbb{R}^n$ \cite{GLSW}.
Furthermore, $\mathbf{F}$ defined on $\mathbb{R}^n$ is
a $C^1$-extension of the curvature $K$ defined on
the space of conformal factors $\Omega^{\mathbb{E}, \mathcal{T}}_{D}(d')$ for $d'\in \mathcal{D}(d)$.
Then we can extend the Euclidean discrete Laplace operator to be defined on $\mathbb{R}^n$, which is the space of
the conformal factors for the discrete conformal class $\mathcal{D}(d)$.
\begin{definition}\label{Euclid Delaunay Laplacian}
Suppose $(S, V)$  is a marked surface with a PL metric $d$.
For a function $f: V\rightarrow \mathbb{R}$ on the vertices, the Euclidean discrete conformal Laplace operator of $d$ on $(S, V)$ is
defined to be the map
\begin{equation*}
\begin{aligned}
\Delta^\mathbb{E}: \mathbb{R}^V&\longrightarrow \mathbb{R}^V\\
f&\mapsto \Delta^\mathbb{E} f,
\end{aligned}
\end{equation*}
where the value of $\Delta^\mathbb{E} f$ at $v_i$ is
\begin{equation*}
\begin{aligned}
\Delta^\mathbb{E} f_i=\sum_{j\sim i}(-\frac{\partial \mathbf{F}_i}{\partial u_j})(f_j-f_i)=-(\widetilde{L}f)_i,
\end{aligned}
\end{equation*}
where $\widetilde{L}_{ij}=\frac{\partial \mathbf{F}_i}{\partial u_j}$ is
an extension of $L_{ij}=\frac{\partial K_i}{\partial u_j}$ for $u\in\Omega^{\mathbb{E}, \mathcal{T}}_{D}(d')$, $d'\in \mathcal{D}(d)$.
\end{definition}

\begin{remark}
Note that $\mathbf{F}$ is $C^1$-smooth on $\mathbb{R}^n$ and $\Delta^{\mathbb{E}, \mathcal{T}}$ is independent of
the Delaunay triangulations of a PL metric,
the operator $\Delta^\mathbb{E}$ is well-defined on $\mathbb{R}^n$.
Furthermore, $\Delta^\mathbb{E}$ is continuous and piecewise smooth on $\mathbb{R}^n$ as a matrix-valued function (\cite{GLSW}, Lemma 5.1).
\end{remark}

\subsection{Euclidean combinatorial Calabi flow with surgery on surfaces}
We can extend the definition of combinatorial Calabi flow to be defined on a conformal class of $d$ as follows.

\begin{definition}
Suppose $(S, V)$  is a marked surface with a PL metric $d_0$.
The Euclidean combinatorial Calabi flow with surgery defined on the conformal class $\mathcal{D}(d_0)$ is
\begin{equation}\label{Euclidean CCF with surgery}
\begin{aligned}
\left\{
  \begin{array}{ll}
    \frac{du_i}{dt}=\Delta^\mathbb{E} \mathbf{F}_i, & \hbox{ } \\
    u_i(0)=0, & \hbox{ }
  \end{array}
\right.
\end{aligned}
\end{equation}
where $\Delta^\mathbb{E}$ is the Euclidean discrete conformal Laplace operator in Definition $\ref{Euclid Delaunay Laplacian}$
defined for the PL metric $d(t)=\mathbf{A}^{-1}((\{p\}, e^{u(t)}))\in \mathcal{D}(d_0)$.
\end{definition}

By Proposition \ref{energy function Euclidean},
$\sum_{i=1}^nu_i$ is invariant along the Euclidean combinatorial Calabi flow with surgery $(\ref{Euclidean CCF with surgery})$.

We can modify the definition of combinatorial Calabi flow with surgery
to search for PL metrics with prescribed combinatorial curvatures
and prove a generalization of Theorem \ref{main theorem for Calabi flow with surgery intro}.

\begin{definition}
Suppose $(S, V)$ is a closed connected marked surface with a PL metric $d_0$ and
$K^*$ is a function defined on $V$.
The Euclidean combinatorial Calabi flow with surgery for $K^*$ defined on the conformal class $\mathcal{D}(d_0)$ is defined to be
\begin{equation}\label{Calabi flow with surgery with target}
\begin{aligned}
\left\{
  \begin{array}{ll}
    \frac{du_i}{dt}=\Delta^\mathbb{E} (\mathbf{F}-K^*)_i, & \hbox{ } \\
    u_i(0)=0. & \hbox{  }
  \end{array}
\right.
\end{aligned}
\end{equation}
\end{definition}
If $K^*=0$, the flow (\ref{Calabi flow with surgery with target}) is reduced to
the flow  (\ref{Euclidean CCF with surgery}).
The Euclidean combinatorial Calabi flow with surgery for $K^*$ (\ref{Calabi flow with surgery with target}) is also a negative gradient flow.
The flow (\ref{Calabi flow with surgery with target}) can be written as
\begin{equation*}
\begin{aligned}
\frac{du_i}{dt}=\Delta^\mathbb{E} (\mathbf{F}-K^*)_i=-\frac{1}{2}\nabla_{u_i}(||\mathbf{F}-K^*||^2)=-\frac{1}{2}\nabla_{u_i}\overline{\mathcal{C}},
\end{aligned}
\end{equation*}
where $\overline{\mathcal{C}}=||\mathbf{F}-K^*||^2$.
This implies that (\ref{Calabi flow with surgery with target})
is the negative gradient flow of $\overline{\mathcal{C}}$
and $\overline{\mathcal{C}}$ is decreasing along the flow (\ref{Calabi flow with surgery with target}).

We have the following generalization of Theorem \ref{main theorem for Calabi flow with surgery intro}
for the Euclidean combinatorial Calabi flow with surgery for $K^*$ (\ref{Calabi flow with surgery with target}).
\begin{theorem}\label{Euclid main theorem for target}
Suppose $(S, V)$ is a closed connected marked surface and $d_0$ is a PL metric on $(S, V)$.
For any $K^*: V\rightarrow (-\infty, 2\pi)$ with $\sum_{v\in V}K^*(v)=2\pi\chi(S)$,
the combinatorial Calabi flow with surgery for $K^*$ $(\ref{Calabi flow with surgery with target})$
exists for all time and converges exponentially fast
to a unique PL metric $d^*\in \mathcal{D}(d_0)$ with combinatorial curvature $K^*$
after finite number of surgeries.
\end{theorem}

\pf
Note that $\sum_{i=1}^nu_i$ is invariant
along the combinatorial Calabi flow with surgery for $K^*$ (\ref{Calabi flow with surgery with target}), then
the solution $u(t)$ of (\ref{Calabi flow with surgery with target})
lies in the hyperplane $\Sigma_0=\{u\in \mathbb{R}^n|\sum_{i=1}^nu_i=0\}$.
As $K^*$ is a given function on $(S, V)$ with $\sum_{i=1}^nK_i^*=2\pi\chi(S)$,
there exists a unique $u^*\in \Sigma_0$ such that $\mathbf{F}(u^*)=K^*$ by Gu-Luo-Sun-Wu's Theorem \ref{Euclidean discrete uniformization}.
Define the Ricci potential function $W^*$ on $R^n$ as
\begin{equation*}
\begin{aligned}
W^*(u)=W(u)-\int_{u^*}^u\sum_{i=1}^nK^*_idu_i=\int_{u^*}^u\sum_{i=1}^n(\mathbf{F}_i-K^*_i)du_i.
\end{aligned}
\end{equation*}
$W^*$ is well-defined and $C^2$ on $\mathbb{R}^n$ with $\nabla W^*=\mathbf{F}-K^*$ by Proposition \ref{energy function Euclidean}.
Furthermore, $W^*$ is convex on $R^n$ and strictly convex on $\Sigma_0$.
So we have $W^*(u^*)=0$, $\nabla W^*(u^*)=0$ and $W^*(u)\geq 0$ for any $u\in \Sigma_0$.

On the other hand, the following lemma is a well-known fact from analysis. The reader could refer to \cite{GX5} (Lemma 4.6) for a proof.
\begin{lemma}\label{lemma on convex function}
Suppose $f(x)$ is a $C^1$-smooth convex function on $\mathbb{R}^{n-1}$ with $\nabla f(x_0)=0$ for some $x_0\in \mathbb{R}^{n-1}$,
$f(x)$ is strictly convex
in a neighborhood of $x_0$, then $\lim_{x\rightarrow \infty}f(x)=+\infty$.
\end{lemma}
By Lemma \ref{lemma on convex function}, we have $\lim_{u\rightarrow \infty}W^*|_{\Sigma_0}=+\infty$
and $W^*|_{\Sigma_0}$ is a nonnegative proper function on $\Sigma_0$.

Set $\psi(t)=W^*(u(t))$. Then
$$\psi'(t)=\nabla_uW^*\cdot \frac{du}{dt}=-(\mathbf{F}-K^*)^T\cdot \widetilde{L} \cdot (\mathbf{F}-K^*)\leq 0$$
by Proposition \ref{energy function Euclidean}.
As $\sum_{i=1}^nu_i$ is invariant along the flow and $\psi(t)\geq 0$, by the properness of $W^*|_{\Sigma_0}$  on $\Sigma_0$,
we have the solution $u(t)$ lies in a compact subset of $\Sigma_0$,
which implies that the solution of the combinatorial Calabi flow (\ref{Calabi flow with surgery with target}) exists for all time.

By the $C^1$-smoothness of $\mathbf{F}$ and the boundness of the solution $u(t)$ of the flow (\ref{Calabi flow with surgery with target}),
to prove the convergence of
the combinatorial Calabi flow with surgery (\ref{Calabi flow with surgery with target}),
we just need to prove the convergence of the curvature.
As the solution $u(t)$ of the combinatorial Calabi flow (\ref{Calabi flow with surgery with target}) lies in a compact subset of $\Sigma_0$
and $\widetilde{L}$ is strictly positive definite on $\Sigma_0$,
we have the first nonzero eigenvalue $\lambda_1$ of $\widetilde{L}$, which is continuous in $u(t)$,
has a uniform positive lower bound $\sqrt{\frac{\lambda}{2}}>0$
along the Euclidean combinatorial Calabi flow with surgery (\ref{Calabi flow with surgery with target}).
So we have
$$\overline{\mathcal{C}}'(t)=-2(\mathbf{F}-K^*)^T\cdot\widetilde{L}^2\cdot(\mathbf{F}-K^*)\leq -2\lambda_1^2(t)\overline{\mathcal{C}}(t)\leq -\lambda \overline{\mathcal{C}}(t),$$
which implies
$$\overline{\mathcal{C}}(t)\leq \overline{\mathcal{C}}(0)e^{-\lambda t}.$$
So the curvature $\mathbf{F}$ converges exponentially fast to $K^*$, which implies that
$u(t)$ converges exponentially fast to $u^*$ by the $C^1$-smoothness of $\mathbf{F}$.

In \cite{Wu}, Wu proved the following result.
\begin{theorem}\label{finiteness}
Suppose $\mathbb{R}^n=\cup D_i$ is an analytic cell decomposition. $f(x)\in C^1(\mathbb{R}^n)$
is analytic on each cell $D_i$, and has a unique minimum
point where $f$ has positive hession. Then its gradient flow $\gamma(t)$, which satisfies $\gamma'(t)=-\nabla f(\gamma(t))$, intersects the
cell faces $D_i$ finitely many times.
\end{theorem}
Combining Theorem \ref{finiteness}, $\overline{\mathcal{C}}(u^*)=0$, (\ref{negative definiteness of critical point}),
the fact that the combinatorial Calabi flow is a negative gradient flow and the fact that
 $\mathbb{R}^n=\cup_{\mathcal{T}} \Omega^{\mathbb{E}, \mathcal{T}}_{D}(d')$ is an analytic cell decomposition of $\mathbb{R}^n$ \cite{GLSW},
we get the finiteness of the surgeries along the combinatorial Calabi flow with surgery (\ref{Calabi flow with surgery with target}).
\qed

\begin{remark}
Ge-Jiang \cite{GJ0} used another extension method introduced by Bobenko-Pinkall-Springborn \cite{BPS}
to extend the combinatorial curvature to the whole space by extending the inner angle of a triangle by constant.
This method ensures the long time existence of the extended combinatorial Yamabe flow and
was then used to study the combinatorial Ricci flow of inversive distance circle packing metrics \cite{GJ1,GJ2,GJ3,GX5}.
Comparing to Ge-Jiang's method to extend combinatorial curvature flow, the extension of Gu-Luo-Sun-Wu
has the following two advantages.
The first is that Ge-Jiang's method is used to study the curvature for a fixed triangulation,
where there may exist no constant curvature polyhedral metric on the triangulation, while Gu-Luo-Sun-Wu's extension ensures the existence of
constant curvature polyhedral metric.
The second is that the extension method used by Ge-Jiang can be used to give a continuous extension of the combinatorial curvature,
but the derivative of curvature may blow up along the flow, the bound of which is crucial for the long time existence and
convergence of combinatorial Calabi flow. The $C^1$-smoothness of $\mathbf{F}$ is crucial for the proof of Theorem $\ref{Euclid main theorem for target}$.
\end{remark}

\section{Hyperbolic combinatorial Calabi flow}\label{section 4}
In this section, we introduce the combinatorial Calabi flow
for piecewise hyperbolic metrics. To handle the singularities that may develop along the flow,
we do surgery by flipping again.
Then we prove the long time existence and convergence of the hyperbolic combinatorial Calabi flow with surgery.
As the results are paralleling to the Euclidean case and the proofs are almost the same,
the results in this section are stated without proof.
\subsection{Hyperbolic combinatorial Calabi flow on triangulated surfaces}
We define the combinatorial Calabi flow for PH metrics similar to the PL case.
\begin{definition}
Suppose $d_0$ is a PH metric on a triangulated surface $(S, V, \mathcal{T})$.
The hyperbolic combinatorial Calabi flow on $(S, V, \mathcal{T})$ is defined as
\begin{equation}\label{hyperbolic Calabi flow for triangulated surface}
\begin{aligned}
\left\{
  \begin{array}{ll}
    \frac{du_i}{dt}=\Delta^{\mathbb{H},\mathcal{T}} K_i, & \hbox{ } \\
    u_i(0)=0, & \hbox{ }
  \end{array}
\right.
\end{aligned}
\end{equation}
where $\Delta^{\mathbb{H},\mathcal{T}}$ is the hyperbolic discrete Laplace operator $(\ref{hyper laplace operator})$
defined by the PH metric $u*d$.
\end{definition}

Similar to the PL case, we have the short time existence of the hyperbolic combinatorial Calabi flow and
the combinatorial curvature $K$ evolves according to
\begin{equation*}
\frac{dK}{dt}=-(\Delta^{\mathbb{H},\mathcal{T}})^2 K
\end{equation*}
along the hyperbolic combinatorial Calabi flow (\ref{hyperbolic Calabi flow for triangulated surface}).
The hyperbolic combinatorial Calabi flow  is also a negative gradient flow of the combinatorial Calabi energy.
\begin{proposition}
The combinatorial Calabi flow $(\ref{hyperbolic Calabi flow for triangulated surface})$
for hyperbolic polyhedral metrics on a triangulated surface $(S, V, \mathcal{T})$
is the negative gradient flow of combinatorial Calabi energy $\mathcal{C}$
and the combinatorial Calabi energy $\mathcal{C}$ is decreasing along the flow.
\end{proposition}

\begin{theorem}\label{stability for triangulated hyper CCF}
If the solution of hyperbolic combinatorial Calabi flow $(\ref{hyperbolic Calabi flow for triangulated surface})$
on a triangulated surface $(S, V, \mathcal{T})$
converges as time tends to infinity, then the limit metric is a hyperbolic metric.
Furthermore, suppose there exists a hyperbolic metric $d^*=u^**d_0$ on $(S, V, \mathcal{T})$,
there exists a constant $\delta>0$
such that if the initial modified Calabi energy $||K(u(0))-K(u^*)||^2<\delta$,
then the hyperbolic combinatorial Calabi flow $(\ref{hyperbolic Calabi flow for triangulated surface})$
on $(S, V, \mathcal{T})$ exists for all time and converges exponentially fast to $u^*$.
\end{theorem}

\subsection{Gu-Guo-Luo-Sun-Wu's work on discrete uniformization theorem}
\begin{definition}[\cite{GGLSW}, Definition 1]\label{hyperbolic discrete conformal}
Two PH metrics $d$, $d'$ on a closed marked surface
$(S, V)$ are discrete conformal if there exists sequences of PH
metrics $d_1=d$, $d_2, \cdots, d_m=d'$ on $(S, V)$ and triangulations $\mathcal{T}_1, \cdots, \mathcal{T}_m$
of $(S, V)$ satisfying
\begin{description}
  \item[(a)] each $\mathcal{T}_i$ is Delaunay in $d_i$,
  \item[(b)] if $\mathcal{T}_i=\mathcal{T}_{i+1}$, there exists a function $u:V\rightarrow \mathbb{R}$,
  called a conformal factor, so that if $e$ is an edge in $\mathcal{T}_i$ with end points $v$ and $v'$,
  then the lengths $x_{d_i}(e)$ and $x_{d_{i+1}}(e)$ of $e$ in metrics $d_i$ and $d_{i+1}$ are related by
  $$\sinh \frac{x_{d_{i+1}}(e)}{2}=e^{u(v)+u(v')}\sinh \frac{x_{d_{i}}(e)}{2},$$
  \item[(c)] if $\mathcal{T}_i\neq \mathcal{T}_{i+1}$, then $(S, d_i)$ is isometric to $(S, d_{i+1})$ by an
  isometry homotopic to the identity in $(S, V)$.
\end{description}
\end{definition}

The space of PH metrics on $(S, V)$ discrete conformal to $d$ is called the conformal class of $d$ and
denoted by $\mathcal{D}(d)$.

Recall the following result of Delaunay triangulations for PH metrics.
\begin{lemma}[\cite{GGLSW}, Proposition 16]\label{hyperbolic surgery}
Suppose $(S, V)$ is a marked surface with a PH metric $d$.
If $\mathcal{T}$ and $\mathcal{T}'$ are Delaunay triangulations of $d$, then there exists
a sequence of Delaunay triangulations $\mathcal{T}_1=\mathcal{T}, \mathcal{T}_2, \cdots, \mathcal{T}_k=\mathcal{T}'$
so that $\mathcal{T}_{i+1}$ is obtained from $\mathcal{T}_i$ by a diagonal switch.
\end{lemma}

A diagonal switch in Lemma \ref{hyperbolic surgery} is referred to a surgery by flipping in the hyperbolic case.

The following discrete uniformization theorem was established by Gu-Guo-Luo-Sun-Wu \cite{GGLSW}.
\begin{theorem}[\cite{GGLSW}, Theorem 3]\label{hyperbolic uniformization}
Suppose $(S, V)$ is a closed connected surface with marked points and
$d$ is a PH metric on $(S, V)$.
Then for any $K^*: V\rightarrow (-\infty, 2\pi)$ with
$\sum_{v\in V}K^*(v)>2\pi\chi(S)$, there exists a unique
PH metric $d'$ on $(S, V)$ so that $d'$ is discrete conformal
to $d$ and the discrete curvature of $d'$ is $K^*$. Furthermore, the discrete
Yamabe flow with surgery associated to curvature $K^*$ having initial value $d$
 converges to $d'$ linearly fast.
\end{theorem}

Denote the Teichim\"{u}ller space of all PH metrics on $(S, V)$ by $T_{hp}(S, V)$ and decorated
Teichim\"{u}ller space of all equivalence class of decorated hyperbolic metrics on $S-V$ by $T_D(S-V)$.
In the proof of Theorem \ref{hyperbolic uniformization}, they prove the following result.

\begin{theorem}[\cite{GGLSW}, Theorem 22, Corollary 24]\label{main result of GGLSW}
There is a $C^1$-diffeomorphism $\mathbf{A}: T_{hp}(S, V)\rightarrow T_D(S, V)$ between $T_{hp}(S, V)$ and $T_D(S-V)$.
Furthermore, the space $\mathcal{D}(d)\subset T_{hp}(S, V)$ of all equivalence classes of PH metrics
discrete conformal to $d$ is $C^1$-diffeomorphic to $\{p\}\times \mathbb{R}^V_{>0}$ under the diffeomorphism $\mathbf{A}$,
where $p$ is the unique hyperbolic metric on $S-V$ determined by the PH metric $d$ on $(S, V)$.
\end{theorem}

Set $u_i=\ln w_i, 1\leq i\leq n$, for $w=(w_1, w_2, \cdots, w_n)\in \mathbb{R}^n_{>0}$ and
$\mathbf{P} =\{x\in (-\infty, 2\pi)^n|\sum_{i=1}^n x_i>2\pi\chi(S)\}$.
Using the map $\mathbf{A}$, Gu-Guo-Luo-Sun-Wu defined the curvature map
\begin{equation*}
\begin{aligned}
\mathbf{F}:\mathbb{R}^n&\rightarrow \mathbf{P}\\
u&\mapsto K_{\mathbf{A}^{-1}(p, w(u))}
\end{aligned}
\end{equation*}
and proved the following property of $\mathbf{F}$.

\begin{proposition}[\cite{GGLSW}]\label{energy function}
There exists a $C^2$-smooth strictly convex function $W: \mathbb{R}^n\rightarrow \mathbb{R}$
so that its gradient $\nabla W$ is $\mathbf{F}$.
\end{proposition}

Theorem \ref{main result of GGLSW} implies that the union of the admissible spaces $\Omega^{\mathbb{H}, \mathcal{T}}_{D}(d')$
of conformal factors such that $\mathcal{T}$ is Delaunay for $d'\in \mathcal{D}(d)$
is $\mathbb{R}^n$.
Similar to the Euclidean case, $\mathbb{R}^n=\cup_{\mathcal{T}} \Omega^{\mathbb{H}, \mathcal{T}}_{D}(d')$ is an analytic cell decomposition of $\mathbb{R}^n$ \cite{GGLSW}.
Furthermore, $\mathbf{F}$ defined on $\mathbb{R}^n$ is
a $C^1$-extension of the curvature $K$ defined on
the space of conformal factors $\Omega^{\mathbb{H}, \mathcal{T}}_{D}(d')$ for $d'\in \mathcal{D}(d)$.
Then we can extend the hyperbolic discrete Laplace operator to be defined on $R^n$, which is the total space of
the conformal factors for the conformal class $\mathcal{D}(d)$.

\begin{definition}
Suppose $(S, V)$  is a marked surface with a PH metric $d$.
For a function $f: V\rightarrow \mathbb{R}$ on the vertices, the hyperbolic discrete conformal Laplace operator is
defined to be the map
\begin{equation*}
\begin{aligned}
\Delta^\mathbb{H}: \mathbb{R}^V&\longrightarrow \mathbb{R}^V\\
f&\mapsto \Delta^\mathbb{H} f,
\end{aligned}
\end{equation*}
where the value of $\Delta^\mathbb{H} f$ at $v_i$ is
\begin{equation*}
\begin{aligned}
\Delta^\mathbb{H} f_i=\sum_{j\sim i}(-\frac{\partial \mathbf{F}_i}{\partial u_j})(f_j-f_i)=-(\widetilde{L}f)_i,
\end{aligned}
\end{equation*}
where $\widetilde{L}_{ij}=\frac{\partial \mathbf{F}_i}{\partial u_j}$ is
an extension of $\frac{\partial K_i}{\partial u_j}$ for a single Delaunay triangulation $\mathcal{T}$ without surgery.
\end{definition}

\begin{remark}
$\mathbf{F}$ is $C^1$-smooth on $\mathbb{R}^n$ and $\Delta^{\mathbb{H}, \mathcal{T}}$ is independent of
the Delaunay triangulations of a PH metric \cite{GGLSW},
the operator $\Delta^\mathbb{H}$ is well-defined on $\mathbb{R}^n$.
Furthermore, $\Delta^\mathbb{H}$ is continuous and piecewise smooth on $\mathbb{R}^n$ as a matrix-valued function (\cite{GGLSW}, Lemma $26$).
\end{remark}

\subsection{Hyperbolic combinatorial Calabi flow with surgery on surfaces}
We can extend the definition of hyperbolic combinatorial Calabi flow to be defined on the discrete conformal class $\mathcal{D}(d_0)$ as follows.

\begin{definition}
Suppose $(S, V)$ is a marked surface with a PH metric $d_0$.
The hyperbolic combinatorial Calabi flow with surgery defined on the conformal class $\mathcal{D}(d_0)$ is defined to be
\begin{equation*}\label{hyperbolic CCF with surgery}
\begin{aligned}
\left\{
  \begin{array}{ll}
    \frac{du_i}{dt}=\Delta^\mathbb{H} \mathbf{F}_i, & \hbox{ } \\
    u_i(0)=0. & \hbox{ }
  \end{array}
\right.
\end{aligned}
\end{equation*}
\end{definition}

We can modify the definition of hyperbolic combinatorial Calabi flow with surgery
to search for PH metrics with prescribing combinatorial curvatures
and prove a generalization of Theorem \ref{main theorem for hyper Calabi flow with surgery intro}.

\begin{definition}
Suppose $(S, V)$ is a marked surface with a PH metric $d_0$ and $K^*$ is a function defined on $V$.
The hyperbolic combinatorial Calabi flow with surgery for $K^*$ defined on the discrete conformal class $\mathcal{D}(d_0)$ is defined as
\begin{equation}\label{hyper Calabi flow with surgery with target}
\begin{aligned}
\left\{
  \begin{array}{ll}
    \frac{du_i}{dt}=\Delta^\mathbb{H} (\mathbf{F}-K^*)_i, & \hbox{ } \\
    u_i(0)=0. & \hbox{ }
  \end{array}
\right.
\end{aligned}
\end{equation}
\end{definition}

The hyperbolic combinatorial Calabi flow with surgery for $K^*$ (\ref{hyper Calabi flow with surgery with target})
is also a negative gradient flow.
The flow (\ref{hyper Calabi flow with surgery with target}) could be written in the following form
\begin{equation*}
\begin{aligned}
\frac{du_i}{dt}=\Delta^\mathbb{E} (\mathbf{F}-K^*)_i=-\frac{1}{2}\nabla_{u_i}(||\mathbf{F}-K^*||^2)=-\frac{1}{2}\nabla_{u_i}\overline{\mathcal{C}},
\end{aligned}
\end{equation*}
where $\overline{\mathcal{C}}=||\mathbf{F}-K^*||^2$.
This implies that (\ref{hyper Calabi flow with surgery with target})
is the negative gradient flow of $\overline{\mathcal{C}}$
and $\overline{\mathcal{C}}$ is decreasing along the flow (\ref{hyper Calabi flow with surgery with target}).

We have the following generalization of Theorem \ref{main theorem for hyper Calabi flow with surgery intro}
for the hyperbolic combinatorial Calabi flow with surgery for $K^*$ (\ref{hyper Calabi flow with surgery with target}).
\begin{theorem}
Suppose $(S, V)$ is a closed connected marked surface and $d_0$ is a PH metric on $(S, V)$.
For any $K^*: V\rightarrow (-\infty, 2\pi)$ with $\sum_{v\in V}K^*(v)>2\pi\chi(S)$,
the hyperbolic combinatorial Calabi flow with surgery $(\ref{hyper Calabi flow with surgery with target})$
exists for all time and converges exponentially fast
to a unique PH metric $d^*\in \mathcal{D}(d_0)$ with combinatorial curvature $K^*$ after finite number of surgeries.
\end{theorem}

\section{Remarks and questions}\label{section 5}
We use the discrete uniformization theorems established in \cite{GGLSW, GLSW}
to prove the long time existence and convergence of the combinatorial Calabi flow with surgery in
Theorem \ref{main theorem for Calabi flow with surgery intro}
and Theorem \ref{main theorem for hyper Calabi flow with surgery intro}.
For the combinatorial Yamabe flow with surgery,
the long time existence and convergence were proved similarly \cite{GGLSW, GLSW}.
Similar results were proved for smooth surface Ricci flow \cite{CH1, H2} and smooth surface Calabi flow \cite{CXX, CP}.
However, the smooth Ricci flow and Calabi flow can be used to prove the uniformization theorem on Riemannian surfaces \cite{CSC, CSC2, CLT}.
So it is interesting to ask the following question.

\textbf{Question: }
Can we reprove the discrete uniformization theorems via the combinatorial Yamabe flow or the combinatorial Calabi flow with surgery?

In \cite{X3, X4}, the second author generalizes the definition of $\alpha$-curvatures in \cite{GX3, GX4} for Thurston's circle
packing metrics to polyhedral metrics on surfaces.
Similar properties of $\alpha$-curvatures on polyhedral surfaces to those in \cite{GX3, GX4} are established in \cite{X3, X4}
using combinatorial Yamabe flow and combinatorial Calabi flow with surgery.
Furthermore, based on the discrete uniformization theorems in \cite{GGLSW, GLSW}, we prove a parameterized discrete uniformization
theorem for $\alpha$-curvatures in \cite{X3, X4}.
\\[10pt]

\textbf{Acknowledgements}\\[8pt]
The first author thanks Professor Jian Sun for introduction and guidance on discrete geometry.
Part of this work was done when the first author was visiting the School of Mathematics and Statistics, Wuhan University.
He would like to thank Wuhan University for its hospitality.
The second author also thanks Professor
Kai Zheng for communications on Calabi flow.
The authors thank Dr. Tianqi Wu for helpful communications.
The research of the second author is supported by Hubei Provincial Natural Science Foundation of China under grant no. 2017CFB681,
National Natural Science Foundation of China under grant no. 61772379 and no. 11301402
and Fundamental Research Funds for the Central Universities.

\end{document}